\documentclass[reqno,A4paper]{amsart}

\usepackage{amsmath}
\usepackage{amssymb}
\usepackage{amsthm}
\usepackage{enumerate}
\usepackage{mathrsfs} 
\usepackage{eqlist}
\usepackage{array}
\usepackage[unicode]{hyperref}
\theoremstyle{plain}
\newtheorem{theorem}{Theorem}[section]

\theoremstyle{definition}

\newtheorem{remark}{\textup{Remark}} 
\newtheorem{example}{\textit{Example}} 

\newtheorem{corollary}{Corollary}[section]
\newtheorem{proposition}{Proposition}[section]
\numberwithin{equation}{section}

\def\R{\mathbb{R}}
\def\lv{\lVert}

\begin{document}

\title[Fundamental inequalities for the iterated Fourier-cosine convolution with Gaussian weight ]%
{Fundamental inequalities for the iterated Fourier-cosine convolution with Gaussian weight and its application}
\author[N. T. H. Phuong \and T. Tuan \and L. T. Minh]%
{Nguyen Thi Hong Phuong* \and Trinh Tuan** \and Lai Tien Minh***}
\date{Accepted 28-April-2025 by \textsc{Mathematica Slovaca} }

\thanks{Published Online: 09-August-2025. \url{https://doi.org/10.1515/ms-2025-0065} }
\newcommand{\acr}{\newline\indent}

\address{\llap{*\,}School for Gifted Students,\acr
                   Hanoi National University of Education,\acr
                   136-Xuan Thuy Rd., Hanoi,\acr VIETNAM.
                   \acr
                  Orcid: \url{https://orcid.org/0009-0007-9025-7924}}
\email{phuongnth@hnue.edu.vn}

\address{\llap{**\,}Department of Mathematics, Faculty of Sciences,\acr
                    Electric Power University,\acr
                    235-Hoang Quoc Viet Rd., Hanoi,\acr
                     VIETNAM.\acr
                    Orcid: \url{https://orcid.org/0000-0002-0376-0238}}
\email{tuantrinhpsac@yahoo.com}

\address{\llap{***\,}Department of Mathematics, \acr
	Hanoi Architectural University,\acr
	Km-10, Nguyen Trai str., Hanoi,\acr
	VIETNAM.\acr
	Orcid: \url{https://orcid.org/0000-0003-2656-8246}}
\email{minhlt@hau.edu.vn}


\begin{abstract}
Derived from the results in [Giang et al.: \emph{Convolutions for the Fourier transforms with geometric variables and applications}, Math. Nachr. 283(12) (2010), 1758–1770], in this paper, we devoted to studying the boundedness properties for the Fourier-cosine convolution weighted by a Gaussian function of the form $\gamma =\exp(-\frac{1}{2}y^2)$ via Young's type theorem and Saitoh's type inequality. New norm estimations in the weighted space are obtained, and the application of the corresponding class of convolutions in Fredholm's second kind of integral equation is discussed. The conditions for the solvability of this equation on $L_1$ space are also found, along with the analysis of an illustrative numerical example, which exemplifies that the present object and method solve cases that are not under the conditions of previously known techniques. 
\end{abstract}

\maketitle

\section{Introduction}
\noindent The theory of convolution in integral transforms has consistently been a dynamic and highly explored field for researchers across applied mathematics, engineering, and physics. At its core, the concept of an integral transform revolves around the properties of the kernel function it utilizes, the most famous of which is the Fourier transform and the convolutions associated with this transform \cite{15,Davies2002integral,GlaeskeHJ2006,grafakos2008,hormander2003analysis,Zayed1996handbook}. 
Based on \cite{5}, the classical Fourier-cosine transform of function $f$, denoted by  ($F_c$) and defined by the integral formula, is as follows.
\begin{equation}\label{1.1}
(F_cf)(y):=\sqrt{\frac{2}{\pi}}\int_{\mathbb R_+}f(x)\cos(xy)dx, \,\,y>0.
\end{equation}\noindent  The Fourier-cosine transform coincides with the Fourier transform when \( f(x) \) is an even function. More generally, the even part of the Fourier transform of \( f(x) \) is equal to the even part of the Fourier-cosine transform of \( f(x) \) within the specified region. In \cite{15}, Churchill  studied the classical convolution of two functions, \( f \) and \( g \), for the Fourier-cosine transform, which is defined by the formula
\begin{equation}\label{1.2}
(f \underset{F_c}{*}g)(x):= \frac{1}{\sqrt{2\pi}} \int_{\mathbb R_+} f(y)\big[g(x+y)+g(|x-y|) \big]dy,\ x>0.
\end{equation}
If \( f \) and \( g \) are \( L_1 \)-Lebesgue integrable functions over \( \mathbb{R}_+ \), then \( (f \underset{F_c}{*} g) \) belongs to \( L_1(\mathbb{R}_+) \), and the factorization identity \( (F_c f)(y) (F_c g)(y)=F_c(f \underset{F_c}{*} g)(y) \) holds for all \( y \in \mathbb{R}_+ \) (see \cite{5}). Following \cite{6},  the authors introduced an alternative definition of the Fourier-cosine integral transform in multi-dimensional space \( \mathbb{R}^n \), denoted by \( (T_c) \) defined by formula $\left(T_c f\right)(x)=\frac{1}{(2 \pi)^{n/ 2}} \displaystyle\int_{\mathbb{R}^n} f(y) \cos (x y)  d y$, and the corresponding convolution of this transform, denoted by \( (f \underset{T_c}{\overset{\gamma}{*}} g) \). Here, \( y = (y_1, y_2, \dots, y_n) \in \mathbb{R}^n \) and $y^2=\sum_{i=1}^{n} y_{i}^2$, with  $\gamma(y)$ is a weight function (see \cite{6}, Theorem 2.5, p. 1762). 
The objective of this paper is study the boundedness properties of the iterated Fourier-cosine convolution \eqref{1.5} in the 1-dimensional case $(n=1)$ associated with $ T_c$ transform, when the weight is a single-variable Gaussian function \( \gamma(y) = e^{-\frac{1}{2} y^2} \).  We rewrite the form of the $T_c$ transform  as follows.
\begin{equation}\label{1.4}
(T_cf)(y)=\frac{1}{(2 \pi)^{1\textfractionsolidus 2}}	\int_{\mathbb{R}}f(x)\cos(xy)dx,\ y\in\mathbb{R},
\end{equation}
\noindent and iterated convolution
\begin{equation}\label{1.5}
(f \underset{T_c}{\overset{\gamma}{*}}g)(x) := \frac{1}{8\pi}	\int_{\mathbb{R}^2}f(u)g(v)G(x,u,v)dudv,  \ x\in\mathbb{R},
\end{equation}
where the Gaussian kernel is defined by
\begin{equation}\label{1.6}
G(x,u,v)=e^{-\frac{1}{2}(x+u+v)^2} +e^{-\frac{1}{2}(x+u-v)^2} +e^{-\frac{1}{2}(x-u+v)^2} +e^{-\frac{1}{2}(x-u-v)^2}.
\end{equation}
Moreover, Theorem 2.5 in \cite{6} has affirmed that if $f, g$ are $L_1$-Lebesgue integrable functions over $\R$,  then convolution $(f \underset{T_c}{\overset{\gamma}{*}}g)$ belongs to $ L_1(\mathbb{R})$, where $T_c$  is defined in  \eqref{1.4}. Additionally, an important property proven in \cite{6} is the factorization identity of the \( T_c \)-transform, which asserts that
\begin{equation}\label{1.7}
T_c(f \underset{T_c}{\overset{\gamma}{*}}g)(y) =\gamma(y)(T_cf)(y)(T_cg)(y), \ y\in\mathbb{R}.
\end{equation}

Integral inequalities serve as fundamental tools for analyzing both the qualitative and quantitative properties of integral transforms and differential equation solutions. In particular, convolution inequalities are essential and in fact indispensable, as numerous integral transforms and differential equation solutions are expressed in terms of convolutions \cite{11,Davies2002integral,Hochstadt,22,21,17}. Among the various convolution-type transforms, the Fourier convolution is undoubtedly the most well-known.
A classical result for an upper bound estimation of the Fourier convolution
is 
$\lv f \underset{F}{*} g  \lv_{L_r (\R)} \leq \lv f\lv_{L_p (\R)} \lv g\lv_{L_q (\R)},$ 
with $p, q, r > 1$ are real numbers such that $1\textfractionsolidus p +1\textfractionsolidus q -1=1\textfractionsolidus r$ and $ f \in L_p (\R),\ g\in L_q(\R)$. This result later became more widely known as Young's convolution inequality for Fourier convolution (see \cite{2}, Chapter V, p. 178). Afterward, Adams and Fournier generalized Young's inequality for the Fourier convolution  (see \cite{3}, Theorem 2.24, p. 33) to include a weight function $w(x)$, as follows
$$\bigg|\int_{\R^n} (f\underset{F}{*} g)(x). w(x)dx \bigg| \leq \lv f\lv_{L_p (\R^n)} \lv g\lv_{L_q (\R^n)} \lv w\lv_{L_r (\R^n)},$$
where $p, q, r>1$ satisfy $\frac{1}{p} +\frac{1}{q}+\frac{1}{r}=2$, for all $f \in L_p (\R^n), g\in L_q (\R^n),$ and $w \in L_r (\R^n)$. One limitation of Young's inequalities is that the obtained result does not hold on $L_2$ Hilbert space.
To address this, in \cite{4}, utilizing the general theorem of reproducing kernels \cite{21}, Saitoh derived a new norm inequality for the iterated Fourier convolution in weighted \( L_p(\mathbb{R}, |\rho_j|) \) space, which takes the following form:
$$\big\lv\big( (F_1 \rho_1) \underset{F}{*}(F_2 \rho_2)\big) (\rho_1 \underset{F}{*} \rho_2)^{\frac{1}{p}-1}\big\lv_{L_p (\R)} \leq \prod\limits_{j=1}^{2}\|F_j\|_{L_p(\mathbb{R},\vert\rho_j\vert)},$$ 
where \( \rho_j \) are non-vanishing functions and \( F_j \in L_p(\mathbb{R}, |\rho_j|) \). Here, the norm of  $F_j$ in the weighted space $L_{p}(\mathbb{R}, \rho_j)$ is understood as $
\| F_j  \|_{L_{p}(\mathbb{R}, \rho_j)}=\big\{\int_{\mathbb R}|F_j(x)|^{p} \rho_j(x) dx\big\}^{1/ p}.$ Notice that, the first version of this inequality appeared in 1984 and was introduced by Saitoh \cite{saitoh1984fundamental} for iterated Laplace convolutions in weighted $L_2$ spaces. Later, Cwikel and Kerman extended Saitoh's result (see \cite{Cwikel1996convolution}) to the generalization case on $L_p, (p>1)$ spaces. Besides, the reverse weighted $L_p$-norm inequality for Fourier convolution has also been investigated and applied to inverse heat source problems \cite{20}. Unlike Young’s inequality, Saitoh’s inequality remains valid for case $p = 2$ (refer \cite{saitoh1984fundamental,4}), which is the most notable distinction between the two inequalities. We refer the reader to \cite{18,9,8,7,TuanHienPhuong} for further alternative versions and applications of these inequalities to various integral transforms.

Our first significant contribution is to present an alternative version of the Young and Saitoh-type theorems for the convolution  \eqref{1.5}, leading to a characterization of the boundedness on specific weighted spaces. The key distinction between the two integral transformations, Fourier-cosine $F_c$ (determined in \eqref{1.1}) and $T_c$ (determined in  \eqref{1.4}), is as follows. For any function \( f \) belonging to \( L_1(\mathbb{R}_{+}) \), the inverse transform of \( F_c \), as defined in \cite{16}, is given by the formula:
\[
(F_c^{-1}f)(x) := \sqrt{\frac{2}{\pi}} \int_{\mathbb{R}_+} f(y) \cos(xy) \, dy,
\]
which holds for all \( x > 0 \). In contrast, the inverse transform of \( T_c \) does not exist for any function \( f \in L_1(\mathbb{R}) \). Moreover, the Fourier-cosine transform \( F_c: L_2(\mathbb{R}_+) \leftrightarrow L_2(\mathbb{R}_+) \) is an automorphism (unitary) on \( L_2(\mathbb{R}_+) \) \cite{16}, meaning that \( F_c \times F_c^{-1} \) is the identity operator. However, this property does not hold for the \( T_c \) transform on the \( L_2 \) space. As a result, it is evident that the iterated convolution \eqref{1.5} and the classical convolution \eqref{1.2} are structurally distinct. The second novel aspect of our work lies in the approach to \( L_1 \)-solvability for the Fredholm integral equation of the second kind, utilizing the structure of convolutions associated with the combination of the \( T_c \) transform.

This paper is organized into four sections as follows. In Section \ref{section2}, building upon the techniques from \cite{24,10} and employing H\"older's inequality, we present an alternative version of the Young-type theorem for the convolution \( (f \underset{T_c}{\overset{\gamma}{*}}g) \) and demonstrate that it is a bounded operator in the \( L_{\infty} (\mathbb{R}) \) space. This leads to the inevitable derivation of the general formula for the Young-type convolution inequality in \eqref{1.5}, with a detailed computation of the sharp upper bound for all inequalities. Section \ref{section3} focuses on the presentation of the Saitoh-type theorem and norm inequalities in weighted \( L_p \) spaces. In Section \ref{section4}, we explore several applications of the convolution \( (f \underset{T_c}{\overset{\gamma}{*}}g) \) to the solvability of certain classes of Fredholm integral equations of the second kind. Specifically, by utilizing the results obtained and leveraging Wiener-L\'evy's theorem \cite{12}, we provide the conditions for the solvability of Fredholm's second-kind integral equations involving operator \eqref{1.5}, and derive explicit \( L_1 \)-solutions. Finally, an example is provided and analyzed at the end of the article to illustrate the results and ensure their validity and applicability.

\section{Young-type inequality for convolution \eqref{1.5}}\label{section2}
\begin{theorem}[Young-type theorem for convolution \eqref{1.5}]\label{Young-typetheorem}
	Let $p,q,r$ and $s$  be real numbers in the open interval $(1,\infty)$ such that $	
	\frac{1}{p}+	\frac{1}{q}+	\frac{1}{r}+	\frac{1}{s}=3.$
	For any functions $f\in L_p(\mathbb{R}),g\in L_q(\mathbb{R}),\gamma\in L_r(\mathbb{R})$ and $h\in L_s(\mathbb{R})$,  we obtain the following estimation
	\begin{equation}
	\label{2.2}
	\left|\displaystyle	\int_{\mathbb{R}}(f\underset{T_c}{\overset{\gamma}{*}}g)(x)h(x)dx\right|\leq \frac{1}{2\pi}\left(\sqrt{\frac{2\pi}{r}}\right)^{1\textfractionsolidus r}\left\|f\right\|_{L_p(\mathbb{R})}\left\|g\right\|_{L_q(\mathbb{R})}\left\|h\right\|_{L_s(\mathbb{R})},
	\end{equation}
	where  $\big(f\underset{T_c}{\overset{\gamma}{*}}g\big)$ is defined by \eqref{1.5} and Gaussian weighted $\gamma(y)=e^{-\frac{1}{2}y^2}$.
\end{theorem}
To prove Theorem \ref{Young-typetheorem}, first, we need the following auxiliary Proposition.
\begin{proposition} 
	\textup{i)} With the definition of the Gaussian kernel given by \eqref{1.6}, we have the following relation 
	\begin{equation}\label{1.8}
	\int_{\mathbb{R}}G(x,u,v)dx=	\int_{\mathbb{R}}G(x,u,v)du=	\int_{\mathbb{R}}G(x,u,v)dv=4\sqrt{2\pi}.  
	\end{equation}
	\textup{ii)} For any functions $f, g$ belonging to $L_1 (\R)$ then $\big(f\underset{T_c}{\overset{\gamma}{*}}g\big) \in L_1 (\R)$, and we obtain
	\begin{equation}\label{1.9}
	\begin{aligned}
	&\Vert f\underset{T_c}{\overset{\gamma}{*}} g\Vert_{L_1(\mathbb{R})}\le \frac{1}{\sqrt{2\pi}} \Vert f\Vert_{L_1(\mathbb{R})}\Vert g\Vert_{L_1(\mathbb{R})}.
	\end{aligned}
	\end{equation}	
\end{proposition} 
\begin{proof} 
	\textup{i)} According to the definition \eqref{1.6} of Gaussian function  and having in mind the well-known Gaussian integral formula (see \cite{13}), then 
	\begin{equation}\label{tichphanGauss}
	\int_{\mathbb{R}} e^{px^2+qx} d x=\sqrt{\frac{\pi}{-p}} e^{-\frac{q^2}{4 p}}\ \text{with}\ \big(p \neq 0,\ \textup{Re}(p) \leq 0\big)\ \text{for}\ \big\{q=0,p=-1 \textfractionsolidus 2\big\},
	\end{equation}  we infer the desired conclusion of equality \eqref{1.8}.\\
	\textup{ii)} Since $f, g\in L_1({\mathbb R})$, by using combination \eqref{1.5},  \eqref{1.6}, \eqref{1.8},  and Fubini theorem, we obtain
	\begin{align*}
	\int_{\mathbb{R}}\Big\vert \big(f\underset{T_c}{\overset{\gamma}{*}} g\big)(x) \Big\vert dx&\leq \frac{1}{8\pi} 	\int_{\mathbb{R}}\Big\{ \int_{\mathbb{R}^2} \vert f(u)\vert \vert g(v)\vert \vert G(x,u,v)\vert dudv\Big\}dx\\
	&= \frac{1}{8\pi} \Big(	\int_{\mathbb{R}}\vert f(u))\vert du\Big) \Big(\int_{\mathbb{R}}\vert g(v)\vert dv\Big) \Big(	\int_{\mathbb{R}}\vert G(x,u,v)\vert dx\Big)\\
	&= \frac{1}{\sqrt{2\pi}}\Vert f\Vert_{L_1(\mathbb{R})}\Vert g\Vert_{L_1(\mathbb{R})}<\infty.
	\end{align*}
	This implies that convolution $(f \underset{T_c}{\overset{\gamma}{*}}g)$ belongs to $L_1(\mathbb{R})$ and  the inequality \eqref{1.9} is valid.  
\end{proof}
\begin{proof}\textit{of Theorem \ref{Young-typetheorem}}.
	Assume that $p_1,q_1,r_1$ and $s_1$ are the conjugate exponentials of $p,q,r, s$, respectively.  This means that $\frac{1}{p}+\frac{1}{p_1}=\frac{1}{q}+\frac{1}{q_1}=\frac{1}{s}+\frac{1}{s_1}=\frac{1}{r}+\frac{1}{r_1}=1$ , together with the assumption of the theorem, we get the correlation between exponential numbers
	as follows
	\begin{equation}\label{2.3}	
	\left\{\begin{array}{l}
	\frac{1}{p_1}+\frac{1}{q_1}+\frac{1}{r_1}+\frac{1}{s_1}=1,\\
	p\left(\frac{1}{q_1}+\frac{1}{r_1}+\frac{1}{s_1}\right)=q\left(\frac{1}{p_1}+\frac{1}{r_1}+\frac{1}{s_1}\right)=r\left(\frac{1}{p_1}+\frac{1}{q_1}+\frac{1}{s_1}\right)=s\left(\frac{1}{p_1}+\frac{1}{q_1}+\frac{1}{r_1}\right)=1.
	\end{array}\right.
	\end{equation}
	By the definition of covolution \eqref{1.5} and Gaussian function \eqref{1.6},  we have
	\begin{equation}\label{2.4}
	\begin{aligned}
	&\left|\displaystyle	\int_{\mathbb{R}}(f\underset{T_c}{\overset{\gamma}{*}}g)(x)h(x)dx\right|\leq\displaystyle	\int_{\mathbb{R}}\left|(f\underset{T_c}{\overset{\gamma}{*}}g)(x)h(x)\right|dx\\
	&=\frac{1}{8\pi}\Bigg\{\int_{\mathbb{R}^3}|f(u)||g(v)||h(x)|e^{-\frac{1}{2}(x+u+v)^2}dudvdx
	+\int_{\mathbb{R}^3}|f(u)||g(v)||h(x)|e^{-\frac{1}{2}(x+u-v)^2}dudvdx\\
	&+\int_{\mathbb{R}^3}|f(u)||g(v)||h(x)|e^{-\frac{1}{2}(x-u+v)^2}du dv dx
	+\int_{\mathbb{R}^3}|f(u)||g(v)||h(x)|e^{-\frac{1}{2}(x-u-v)^2}du dv dx\Bigg\}\\
	&=\frac{1}{8\pi}\bigg\{A(u,v,x)+B(u,v,x)+C(u,v,x)+D(u,v,x)\bigg\}.\end{aligned}
	\end{equation}
	For short, we set \begin{align*} 
	A_1(u,v,x)&=|g(v)|^{\frac{q}{p_1}}|h(x)|^{\frac{s}{p_1}}\left(e^{-\frac{1}{2}(x+u+v)^2}\right)^{\frac{r}{p_1}}\in L_{p_1}(\mathbb{R}^3),\\
	A_2(u,v,x)&=|f(u)|^{\frac{p}{q_1}}|h(x)|^{\frac{s}{q_1}}\left(e^{-\frac{1}{2}(x+u+v)^2}\right)^{\frac{r}{q_1}}\in L_{q_1}(\mathbb{R}^3),\\
	A_3(u,v,x)&=|f(u)|^{\frac{p}{s_1}}|g(v)|^{\frac{q}{s_1}}\left(e^{-\frac{1}{2}(x+u+v)^2}\right)^{\frac{r}{s_1}}\in L_{p_1}(\mathbb{R}^3),\\
	A_4(u,v,x)&=|f(u)|^{\frac{p}{r_1}}|g(v)|^{\frac{q}{r_1}}|h(x)|^{\frac{s}{r_1}}\in L_{r_1}(\mathbb{R}^3).
	\end{align*} 
	Due to \eqref{2.3},  we infer that $\prod\limits_{i=1}^{4}A_i(u,v,x)=|f(u)|\,|g(v)|\,|h(x)|\,e^{-\frac{1}{2}(x+u+v)^2}.$	Therefore	
	\begin{equation}\label{2.5}A(u,v,x)=\int_{\mathbb{R}^3}\prod\limits_{i=1}^{4}A_i(u,v,x)dudvdx.\end{equation}
	Moreover, since $\frac{1}{p_1}+\frac{1}{q_1}+\frac{1}{r_1}+\frac{1}{s_1}=1$, applying Hölder's inequality for \eqref{2.5}, we deduce that	
	\begin{equation}\label{2.6}
	\begin{aligned}
	A(u,v,x)\leq&\left\{ \int_{\mathbb{R}^3}|A_1(u,v,x)|^{p_1}dudvdx\right\}^{\frac{1}{p_1}}\left\{ \int_{\mathbb{R}^3}|A_2(u,v,x)|^{q_1}dudvdx\right\}^{\frac{1}{q_1}}\\
	&\left\{ \int_{\mathbb{R}^3}|A_3(u,v,x)|^{r_1}dudvdx\right\}^{\frac{1}{r_1}}\left\{ \int_{\mathbb{R}^3}|A_4(u,v,x)|^{s_1}dudvdx\right\}^{\frac{1}{s_1}}.
	\\ &=\|A_1\|_{L_{p_1}(\mathbb{R}^3)}\|A_2\|_{L_{q_1}(\mathbb{R}^3)}\|A_3\|_{L_{r_1}(\mathbb{R}^3)}\|A_4\|_{L_{s_1}(\mathbb{R}^3)}.
	\end{aligned}\end{equation}
	Based on the assumption of $g\in L_q(\mathbb{R})$, $e^{-\frac{1}{2}(x+u+v)^2}\in L_r(\mathbb{R})$ and $h\in L_s(\mathbb{R})$,  using Fubini’s
	theorem, we obtain $$\|A_1\|_{L_{p_1}(\mathbb{R}^3)}^{p_1}=\left(	\int_{\mathbb{R}}|g(v)|^q dv\right)\left(	\int_{\mathbb{R}}|h(x)|^s dx\right)\left[	\int_{\mathbb{R}}\left(e^{-\frac{1}{2}(x+u+v)^2}\right)^r du\right].$$ Using \eqref{tichphanGauss} for the case $p=\frac{-1}{2}r$, $q=0$ together with the change of variables theorem, we deduce that $$\int_{\mathbb{R}} e^{-\frac{1}{2}r(x+u+v)^2}du=\sqrt{\frac{2\pi}{r}}.$$  So the above equation can be simplified to
	$\left\|A_1 \right\|^{p_1}_{L_{p_1}(\mathbb{R}^3)}=\sqrt{\frac{2\pi}{r}}\left\|g\right\|^q_{L_q(\mathbb{R})}\left\|h\right\|^s_{L_s(\mathbb{R})}.$
	Therefore, we obtain ${L_{p_1}(\mathbb{R}^3)}$-norm estimation for the operator $A_1$ as follows
	\begin{equation}
	\label{2.7}
	\left\|A_1 \right\|_{L_{p_1}(\mathbb{R}^3)}=\left(\sqrt{\frac{2\pi}{r}}\right)^{\frac{1}{p_1}}\left\|g\right\|^{\frac{q}{p_1}}_{L_q(\mathbb{R})}\left\|h\right\|^{\frac{s}{p_1}}_{L_s(\mathbb{R})}.
	\end{equation}
	Similar to what we did with the evaluation \eqref{2.7} of $A_1$, we also get the norm estimation of $A_2,A_3, A_4$
	\begin{equation}
	\label{2.8}
	\left\|A_2 \right\|_{L_{q_1}(\mathbb{R}^3)}=\left(\sqrt{\frac{2\pi}{r}}\right)^{\frac{1}{q_1}}\left\|f\right\|^{\frac{p}{q_1}}_{L_p(\mathbb{R})}\left\|h\right\|^{\frac{s}{q_1}}_{L_s(\mathbb{R})},
	\end{equation}
	\begin{equation}
	\label{2.9}
	\left\|A_3 \right\|_{L_{s_1}(\mathbb{R}^3)}=\left(\sqrt{\frac{2\pi}{r}}\right)^{\frac{1}{s_1}}\left\|f\right\|^{\frac{p}{s_1}}_{L_p(\mathbb{R})}\left\|g\right\|^{\frac{q}{s_1}}_{L_q(\mathbb{R})},
	\end{equation}
	\begin{equation}
	\label{2.10}
	\left\|A_4 \right\|_{L_{r_1}(\mathbb{R}^3)}=\left(\sqrt{\frac{r}{2\pi}}\right)^{\frac{1}{r_1}}\left(\sqrt{\frac{2\pi}{r}}\right)^{\frac{1}{r_1}}\left\|f\right\|^{\frac{p}{r_1}}_{L_p(\mathbb{R})}\left\|g\right\|^{\frac{q}{r_1}}_{L_q(\mathbb{R})}\left\|h\right\|^{\frac{s}{r_1}}_{L_s(\mathbb{R})}.
	\end{equation}
	Coupling \eqref{2.7}, \eqref{2.8}, and \eqref{2.9}, \eqref{2.10}, together with \eqref{2.3}, we obtain
	\begin{equation*}
	\left\|A_1 \right\|_{L_{p_1}(\mathbb{R}^3)}\left\|A_2 \right\|_{L_{q_1}(\mathbb{R}^3)}\left\|A_3 \right\|_{L_{s_1}(\mathbb{R}^3)}\left\|A_4 \right\|_{L_{r_1}(\mathbb{R}^3)}\leq\sqrt{\frac{2\pi}{r}}\left(\sqrt{\frac{r}{2\pi}}\right)^{\frac{1}{r_1}}\left\|f\right\|_{L_p(\mathbb{R})}\left\|g\right\|_{L_q(\mathbb{R})}\left\|h\right\|_{L_s(\mathbb{R})}.
	\end{equation*}
	Combining \eqref{2.6}, we implies that 
	\begin{equation}\label{2.11}
	A(u,v,x)\leq \sqrt{\frac{2\pi}{r}}\left(\sqrt{\frac{r}{2\pi}}\right)^{\frac{1}{r_1}}\|f\|_{L_p(\mathbb{R}^3)}\|g\|_{L_q(\mathbb{R}^3)}\|h\|_{L_s(\mathbb{R}^3)}.
	\end{equation}
	By repeating the above arguments, we also obtain similar evaluations for $B(u,v,x); C(u,v,x)$ and $D(u,v,x)$ as follows
	\begin{equation}\label{2.12}
	\begin{aligned}
	B(u,v,x)&\leq \sqrt{\frac{2\pi}{r}}\left(\sqrt{\frac{r}{2\pi}}\right)^{\frac{1}{r_1}}\|f\|_{L_p(\mathbb{R}^3)}\|g\|_{L_q(\mathbb{R}^3)}\|h\|_{L_s(\mathbb{R}^3)},\\
	C(u,v,x)&\leq \sqrt{\frac{2\pi}{r}}\left(\sqrt{\frac{r}{2\pi}}\right)^{\frac{1}{r_1}}\|f\|_{L_p(\mathbb{R}^3)}\|g\|_{L_q(\mathbb{R}^3)}\|h\|_{L_s(\mathbb{R}^3)},\\
	D(u,v,x)&\leq \sqrt{\frac{2\pi}{r}}\left(\sqrt{\frac{r}{2\pi}}\right)^{\frac{1}{r_1}}\|f\|_{L_p(\mathbb{R}^3)}\|g\|_{L_q(\mathbb{R}^3)}\|h\|_{L_s(\mathbb{R}^3)}.\\
	\end{aligned}
	\end{equation}
	Finally, combining \eqref{2.4},\eqref{2.11}, and \eqref{2.12}, we deduce the estimation as in the conclusion of the theorem. 
\end{proof}

\begin{corollary}[Young-type inequality for convolution \eqref{1.5}]\label{Young-inequality}
	If $p,q,r$ and $s\geq 1$ be real numbers, satisfying $\frac{1}{p}+\frac{1}{q}+\frac{1}{r}=2+\frac{1}{s}$.  If $f\in L_p(\mathbb{R}), g\in L_q(\mathbb{R})$ and $\gamma\in L_r(\mathbb{R})$  then the covolution
	$(f\underset{T_c}{\overset{\gamma}{*}}g)$ belongs to $L_s(\mathbb{R})$.  Furthermore,  we obtain the following inequality 
	\begin{equation}
	\|f\underset{T_c}{\overset{\gamma}{*}}g\|_{L_s(\mathbb{R})}\leq \frac{1}{2\pi}\left(\sqrt{\frac{2\pi}{r}}\right)^{1\textfractionsolidus r}\|f\|_{L_p(\mathbb{R})}\|g\|_{L_q(\mathbb{R})}.
	\label{2.13}
	\end{equation}
\end{corollary}
\begin{proof}
	With $p=q=r=s=1$ has already been proved in norm estimation \eqref{1.9}.  Therefore,  we need to show this corollary holds true for the case where $p,q,r$,  and $s$ are greater than 1. Indeed,  let $s_1$ be the conjugate exponent of $s$, i.e.  $\frac{1}{s}+\frac{1}{s_1}=1$. From the assumptions of Corollary \ref{Young-inequality}, we have $\frac{1}{p}+\frac{1}{q}+\frac{1}{r}+\frac{1}{s_1}=3$, which shows the numbers $p,q,r$, and $s_1$ satisfy the conditions of Theorem \ref{Young-typetheorem} (with $s$ being replaced
	by $s_1$). Choosing the function $h(x)=\big(f \underset{T_c}{\overset{\gamma}{*}}g\big)^\alpha (x),$ where $\alpha$ is a constant that only depends on $s_1$ such that $h(x)\in L_{s_1}(\mathbb{R})$, (constant $\alpha$  certainly exists and we will show how to choose it at the end of this proof). Obviously $\big(f\underset{T_c}{\overset{\gamma}{*}}g\big)(x)\in L_{\alpha s_1}(\mathbb{R})$ and we have
	$$\begin{aligned}
	\|h\|_{L_{s+1}(\mathbb{R})}&=\left\{	\int_{\mathbb{R}}|h(x)|^{s_1}dx\right\}^{\frac{1}{s_1}}=\left\{	\int_{\mathbb{R}}\bigg|(f\underset{T_c}{\overset{\gamma}{*}}g)^\alpha(x)\bigg|^{s_1}dx\right\}^{\frac{1}{s_1}}=
	\left\{	\int_{\mathbb{R}}\bigg|(f\underset{T_c}{\overset{\gamma}{*}}g)^\alpha(x)\bigg|^{\alpha s_1}dx\right\}^{\frac{1}{\alpha s_1}\cdot \alpha}\\&=\|f\underset{T_c}{\overset{\gamma}{*}}g \|^\alpha_{L_{\alpha s_1}(\mathbb{R})}.
	\end{aligned}$$
	By applying Theorem \ref{Young-typetheorem} for $\|h(x)\|_{L_{s_1}(\mathbb{R})}=\|f\underset{T_c}{\overset{\gamma}{*}}g\|_{L_{\alpha s_1}(\mathbb{R})}^{\alpha}$ with $s$ being replaced
	by $s_1$, together with estimate \eqref{2.2}, we deduce that
	$$\left|	\int_{\mathbb{R}}\big(f\underset{T_c}{\overset{\gamma}{*}}g\big)^{\alpha+1}(x)\mathrm{d}x\right|\leq \frac{1}{2\pi}\left(\sqrt{\frac{2\pi}{r}}\right)^{1\textfractionsolidus r}\|f\|_{L_p(\mathbb{R})}\|g\|_{L_q(\mathbb{R})}\|
	f\underset{T_c}{\overset{\gamma}{*}}g
	\|_{L_{\alpha s_1}(\mathbb{R})}^{\alpha}.$$ This is equivalent to
	\begin{equation}\label{2.14}
	\left|	\int_{\mathbb{R}}\big(f\underset{T_c}{\overset{\gamma}{*}}g\big)^{\alpha s_1\cdot \frac{\alpha+1}{\alpha s_1}}(x)dx\right|\leq \frac{1}{2\pi}\left(\sqrt{\frac{2\pi}{r}}\right)^{1\textfractionsolidus r}\|g\|_{L_q(\mathbb{R})}\|
	f\underset{T_c}{\overset{\gamma}{*}}g
	\|_{L_{\alpha s_1}(\mathbb{R})}^{\alpha}.	
	\end{equation}
	Now, we just need to choose $\alpha$ to satisfy $\frac{\alpha+1}{\alpha s_1}=1$. This means that $\alpha=\frac{1}{s_1 -1}$ and \eqref{2.14} becomes the following
	\begin{equation}\label{2.15}
	\|
	f\underset{T_c}{\overset{\gamma}{*}}g
	\|_{L_{\alpha s_1}(\mathbb{R})}^{\alpha(s_1-1)}\leq \frac{1}{2\pi}\left(\sqrt{\frac{2\pi}{r}}\right)^{1\textfractionsolidus r}\|f\|_{L_p(\mathbb{R})}\|g\|_{L_q(\mathbb{R})}.
	\end{equation}
	Since $1\textfractionsolidus s +1\textfractionsolidus s_1 =1$ i.e. $s_1 = \frac{s}{s-1}$.  Therefore $\alpha=\frac{1}{\frac{s}{s-1}-1}=s-1$ implies that $\alpha s_1=(s-1)\frac{s}{s-1}=s$, which is enough to show that $L_s(\mathbb{R})\equiv L_{\alpha s_1}(\mathbb{R})$ and $\alpha(s_1-1)=1$. Combining with inequality \eqref{2.15},  we arrive at the conclusion of this corollary.
\end{proof} Another way to prove Corollary \ref{Young-inequality} is to skillfully 
apply Riesz's representation theorem  \cite{2} skillfully to the bounded linear functional on dual spaces (\cite{14}, proof of Corollary 4.2, page 1691). In addition, the following results are straightforward to get from the condition $\frac{1}{p}+\frac{1}{q}+\frac{1}{r}=2+\frac{1}{s},$ with  $p,q,r,s\geq 1$ and estimate \eqref{2.13}.
\begin{remark}
	If $s=p>1$ and $q=r=1$,  then the following norm estimate holds
	\begin{equation}\label{2.15a}
	\|
	f\underset{T_c}{\overset{\gamma}{*}}g
	\|_{L_{p}(\mathbb{R})}\leq \frac{1}{\sqrt{2\pi}}\|f\|_{L_p(\mathbb{R})}\|g\|_{L_1(\mathbb{R})}.
	\end{equation}
	If $s=r>1$ and $p=q=1$,  then we deduce the following norm estimate
	\begin{equation}\label{2.15b}
	\|
	f\underset{T_c}{\overset{\gamma}{*}}g
	\|_{L_{r}(\mathbb{R})}\leq \frac{1}{2\pi}\left(\sqrt{\frac{2\pi}{r}}\right)^{1\textfractionsolidus r}\|f\|_{L_1(\mathbb{R})}\|g\|_{L_1(\mathbb{R})}.
	\end{equation}	
\end{remark}
\noindent What about the case $s=\infty$? We consider the boundedness of operator \eqref{1.5} in $L_{\infty}(\R)$ via the following theorem.
\begin{theorem}\label{Định lý 2.2}
	Suppose that  $p,q,r>1$ satisfy $\frac{1}{p}+\frac{1}{q}+\frac{1}{r}=2$. For any functions $f\in L_{p}(\mathbb{R}), g\in L_q(\mathbb{R}),$ and $\gamma(y)=e^{-\frac{1}{2}y^2}\in L_r(\mathbb{R})$. Then convolution $\eqref{1.5}$ is a bounded operator $\forall x\in \mathbb{R}$.  Furthermore, the following inequality holds
	\begin{equation}\label{2.16}
	\|
	f\underset{T_c}{\overset{\gamma}{*}}g
	\|_{L_{\infty}(\mathbb{R})}\leq \frac{1}{2\pi}\left(\sqrt{\frac{2\pi}{r}}\right)^{1\textfractionsolidus r}\|f\|_{L_p(\mathbb{R})}\|g\|_{L_q(\mathbb{R})}.	
	\end{equation}
	Here, the norm of functions in $L$-infinity is understood as $\|f\|_{L_{\infty}(\mathbb{R})}:=\mathrm{ess}\sup\limits_{x\in\mathbb{R}}|f(x)|$.
\end{theorem}
\begin{proof}
	Due to \eqref{1.5}, \eqref{1.6} and \eqref{1.8} allow us to recognize \begin{align*}
	&\bigg|\big(f\underset{T_c}{\overset{\gamma}{*}}g
	\big)(x)\bigg|= \frac{1}{8\pi}\left|\int_{\mathbb{R}^2}f(u)g(v)G(x,u,v)du dv\right|\leq \frac{1}{8\pi}\displaystyle\int_{\mathbb{R}^2}\left|f(u)\right|\left|g(v)\right|\left|G(x,u,v)\right|du dv\\
	&\leq \frac{1}{8\pi}\bigg\{\displaystyle\int_{\mathbb{R}^2}\left|f(u)\right|\left|g(v)\right|e^{-\frac{1}{2}(x+u+v)^2}du dv+
	\int_{\mathbb{R}^2}\left|f(u)\right|\left|g(v)\right|e^{-\frac{1}{2}(x+u-v)^2}du dv
	\\
	&+\int_{\mathbb{R}^2}\left|f(u)\right|\left|g(v)\right|e^{-\frac{1}{2}(x-u+v)^2}\mathrm{d}u\mathrm{d}v+\int_{\mathbb{R}^2}\left|f(u)\right|\left|g(v)\right|e^{-\frac{1}{2}(x-u-v)^2}du dv
	\bigg\}=\frac{1}{8\pi}\{I_1 + I_2 + I_3 + I_4 \},
	\end{align*}
	In what follows,  the integrals $I_i, i=\overline{1,4}$ will be estimated.  We first set functions $T_i(x,u,v),$ with $i=\overline{1,4 }$ as
	\begin{align*}
	T_1(x,u,v)&=|g(v)|^{\frac{q}{p_1}}\big\{ e^{-\frac{1}{2}(x+u+v)^2} \big\}^{\frac{r}{p_1}}\in L_{p_1}(\mathbb{R}^2),x\in \mathbb{R}\\
	T_2(x,u,v)&=|f(u)|^{\frac{p}{q_1}}\big\{ e^{-\frac{1}{2}(x+u+v)^2} \big\}^{\frac{r}{q_1}}\in L_{q_1}(\mathbb{R}^2),x\in \mathbb{R}\\
	T_3(x,u,v)&=|f(u)|^{\frac{p}{r_1}}|g(v)|^{\frac{q}{r_1}}
	\in L_{r_1}(\mathbb{R}^2),x\in \mathbb{R},
	\end{align*}
	where $\frac{1}{p}+\frac{1}{p_1}=1$, $\frac{1}{q}+\frac{1}{q_1}=1$ and $\frac{1}{r}+\frac{1}{r_1}=1$. On the other hand,  owing to the condition $\frac{1}{p}+\frac{1}{q}+\frac{1}{r}=2$, we obtain
	\begin{equation}\label{2.17}
	\left\{\begin{array}{l}
	\frac{1}{p_1}+\frac{1}{q_1}+\frac{1}{r_1}=1,\\
	p\left(\frac{1}{q_1}+\frac{1}{r_1}\right)=q\left(\frac{1}{p_1}+\frac{1}{r_1}\right)=r\left(\frac{1}{p_1}+\frac{1}{q_1}\right)=1.
	\end{array}\right.
	\end{equation}	
	Under the condition	\eqref{2.17}, we have $T_1(x,u,v)T_2(x,u,v)T_3(x,u,v)=|f(u)||g(v)|e^{-\frac{1}{2}(x+u+v)^2}.$ This indicates that
	\begin{equation}\label{2.18}
	I_1=\int_{\mathbb{R}^2} T_1(x,u,v)T_2(x,u,v)T_3(x,u,v) dudv,x>0.
	\end{equation}
	Since $\frac{1}{p_1}+\frac{1}{q_1}+\frac{1}{r_1}=1$, by applying H\"older inequality for the right-hand side of  \eqref{2.18}, we obtain
	\begin{equation}\label{2.19}
	I_1 \leq \| T_1\|_{L_{p_1}(\mathbb{R}^2)}\|T_2\| _{L_{q_1}(\mathbb{R}^2)}\|T_3\|_{L_{r_1}(\mathbb{R}^2)}.
	\end{equation}
	By proceeding similarly to the proof of Theorem \ref{Young-typetheorem}, we get the norm estimation  as follows$$\begin{aligned}
	\| T_1\|_{L_{p_1}(\mathbb{R}^2)}&\leq \left(\sqrt{\frac{2\pi}{r}}\right)^{\frac{1}{p_1}}\|g\|_{L_q(\mathbb{R})}^{\frac{q}{p_1}}\\
	\| T_2\|_{L_{q_1}(\mathbb{R}^2)}&\leq \left(\sqrt{\frac{2\pi}{r}}\right)^{\frac{1}{q_1}}\|f\|_{L_p(\mathbb{R})}^{\frac{p}{q_1}}\\
	\| T_3\|_{L_{r_1}(\mathbb{R}^2)}&= \left(\sqrt{\frac{r}{2\pi}}\right)^{\frac{1}{r_1}} \left(\sqrt{\frac{2\pi}{r}}\right)^{\frac{1}{r_1}}\|f\|_{L_p(\mathbb{R})}^{\frac{p}{r_1}}\|g\|_{L_q(\mathbb{R})}^{\frac{q}{r_1}}.
	\end{aligned}$$
	Coupling \eqref{2.17},\eqref{2.19}, we infer that $$I_1\leq \|T_1\|_{L_{p_1}(\mathbb{R}^2)}\|T_2\|_{L_{q_1}(\mathbb{R}^2)}
	\|T_3\|_{L_{r_1}(\mathbb{R}^2)}	
	\leq \left(\sqrt{\frac{r}{2\pi}}\right)^{\frac{1}{r_1}}\sqrt{\frac{2\pi}{r}}\|f\|_{L_p(\mathbb{R})}\|g\|_{L_q(\mathbb{R})}.$$
	Again, same as done with $I_1$, we obtain $I_i(x)\leq \left(\sqrt{\frac{r}{2\pi}}\right)^{\frac{1}{r_1}} \sqrt{\frac{2\pi}{r}}\|f\|_{L_p(\mathbb{R})}\|g\|_{L_q(\mathbb{R})},$ with $i=\overline{2,4}$, combining with \eqref{1.8},  it is easy to verify $$\big|(f\underset{T_c}{\overset{\gamma}{*}}g)(x)\big|\leq \frac{1}{8\pi}\sum\limits_{i=1}^{4} I_i\leq \frac{1}{8\pi} 4\sqrt{\frac{2\pi}{r}} \left(\sqrt{\frac{r}{2\pi}}\right)^{1-1\textfractionsolidus r}	\|f\|_{L_p(\mathbb{R})}\|g\|_{L_q(\mathbb{R})}<\infty$$  for any $f\in L_p(\mathbb{R})$ and $g\in L_q(\mathbb{R})$. This implies that $\big(f\underset{T_c}{\overset{\gamma}{*}}g\big)(x)$ is a bounded operator $\forall x\in \mathbb{R}$ and $\mathrm{ess}\sup\limits_{x\in\mathbb{R}}\big|(f\underset{T_c}{\overset{\gamma}{*}}g)(x)\big|<\infty, $  and infers
	the desired conclusion of inequality \eqref{2.16}.
\end{proof}
\noindent From the inequalities \eqref{1.9}, \eqref{2.13}, and \eqref{2.16}, we deduce that the operator defined in \eqref{1.5} induces a bilinear mapping
$
\underset{T_c}{\overset{\gamma}{*}} : L_p(\mathbb{R}) \times L_q(\mathbb{R}) \to L_s(\mathbb{R}),
$
such that for all \( f \in L_p(\mathbb{R}) \) and \( g \in L_q(\mathbb{R}) \), the mapping \((f, g) \mapsto (f \underset{T_c}{\overset{\gamma}{*}} g)\) is separately continuous in each variable. Furthermore, the operator \eqref{1.5} is well-defined and bounded in the space \(L_s(\mathbb{R})\) for every \(s \in [1, \infty]\).
On the other hand, it is evident that
$
\left( \sqrt{\frac{2\pi}{r}} \right)^{1/r} \leq \sqrt{2\pi} \exp\left( \frac{1}{2e} \right) = \text{Const}.
$
Consequently, by inequality \eqref{2.16}, we obtain the estimate
\[
\| f \underset{T_c}{\overset{\gamma}{*}} g \|_{L_{\infty}(\mathbb{R})} \leq \frac{\exp\left( \frac{1}{2e} \right)}{\sqrt{2\pi}} \|f\|_{L_p(\mathbb{R})} \|g\|_{L_q(\mathbb{R})},
\]
which shows that the convolution operator \eqref{1.5} is uniformly bounded on \(L_\infty(\mathbb{R})\) with respect to each function \( f \in L_p(\mathbb{R}) \) and \( g \in L_q(\mathbb{R}) \).
However, it is important to observe that the estimates \eqref{2.2}, \eqref{2.13}, and \eqref{2.16} are no longer valid in the Hilbert space \( L_2(\mathbb{R}) \); that is, in the specific case when \( p = q = r = s = 2 \).

\section{Saitoh-type inequality for convolution \eqref{1.5}}\label{section3}
\noindent Throughout in this section,  we shall make frequent use of weighted Lebesgue spaces $L_p(\mathbb{R}_+, \rho(x)),$ with  $1\leq p\leq \infty$ with respect to a positive measure $\rho(x)dx,$ equipped with the norm $$\|f\|_{L_p(\mathbb{R}_+,\rho)}=\bigg\{\int_{\R_+} |f(x)|^p \rho(x) dx\bigg\}^{1\textfractionsolidus p}<\infty.$$
\begin{theorem}[Saitoh-type inequality for convolution \eqref{1.5}] \label{Định lý 3.1} Assume that $\rho_i,$ with $i=\{1,2\}$ are non-vanishing positive functions such that  $\big(\rho_1\underset{T_c}{\overset{\gamma}{*}}\cdot\rho_2\big)$ is well-defined for any functions $F_i\in L_p(\mathbb{R},\rho_i),p>1$. Then, we have the following estimate
	\begin{equation}	\label{3.1}
	\left\|\left( \left( F_1\rho_1 \right)\underset{T_c}{\overset{\gamma}{*}}\left( F_2\rho_2 \right)\right) (x)\cdot\left( \rho_1\underset{T_c}{\overset{\gamma}{*}}\rho_2 \right)^{\frac{1}{p}-1}(x)\right\|_{L_p(\mathbb{R})}\leq \left(\frac{1}{\sqrt{2\pi}}\right)^{1\textfractionsolidus p}\prod\limits_{i=1}^{2}\|F_i\|_{L_p(\mathbb{R},\rho_i)},
	\end{equation}
	where convolution $\big(f\underset{T_c}{\overset{\gamma}{*}}g\big)$ is defined by \eqref{1.5}.
\end{theorem}
\begin{proof}
	Putting $\mathcal{M}^p(x)=\left\|\left( \left( F_1\rho_1 \right)\underset{T_c}{\overset{\gamma}{*}}\left( F_2\rho_2 \right)\right) (x)\cdot( \rho_1\underset{T_c}{\overset{\gamma}{*}}\rho_2 )^{\frac{1}{p}-1}(x)\right\|_{L_p(\mathbb{R})}^p.$
	Based on $\eqref{1.5}$,  we obtain	
	\begin{equation}\label{3.2}
	\begin{aligned}	
	&\mathcal{M}^p(x)=\int_{\mathbb{R}}\bigg|\bigg( \big( F_1\rho_1 \big)\underset{T_c}{\overset{\gamma}{*}}\big( F_2\rho_2 \big)\bigg) (x)\cdot( \rho_1\underset{T_c}{\overset{\gamma}{*}}\rho_2 )^{\frac{1}{p}-1}(x)\bigg|^p dx\\
	&=\left(\frac{1}{8\pi}\right)^p\left(\frac{1}{8\pi}\right)^{1-p}\int_{\mathbb{R}} \bigg\{
	\bigg|\int_{\mathbb{R}^2} (F_1\rho_1)(u)(F_2\rho_2)(v)G(x,u,v) dudv\bigg|^p\times
	\bigg|\int_{\mathbb{R}^2} \rho_1(u)\rho_2(v)G(x,u,v) dudv\bigg|^{1-p}
	\bigg\}\\
	&\leq \frac{1}{8\pi} \int_{\mathbb{R}}\bigg\{\bigg(
	\int_{\mathbb{R}^2} \big|(F_1\rho_1)(u)\big|\big|(F_2\rho_2)(v)\big|G(x,u,v)dudv\bigg)^p\times
	\bigg(
	\int_{\mathbb{R}^2}\rho_1(u)\rho_2(v)G(x,u,v)dudv
	\bigg)^{1-p}
	\bigg\}dx.
	\end{aligned}
	\end{equation}
	Using H\"older’s inequality for the pair of conjugate exponents $p,q,$ i.e. $(1\textfractionsolidus p + 1\textfractionsolidus q =1),$ with $p>1$   implies that
	\begin{equation}\label{3.3}
	\begin{aligned}
	&\int_{\mathbb{R}}\big|\big(F_1\rho_1\big)(w)\big|\big|\big(F_2\rho_2\big)(w)\big|G(x,u,v)dudv\\
	&\leq \bigg\{\int_{\mathbb{R}^2}|F_1(u)|^p\rho_1(u)|F_2(v)|^p\rho_2(v)G(x,u,v)\mathrm{d}u\mathrm{d}v\bigg\}^{\frac{1}{p}}\times
	\bigg\{
	\int_{\mathbb{R}^2} \big|(F_1\rho_1)(u)\big|\big|(F_2\rho_2)(v)\big|G(x,u,v)dudv\bigg\}^{\frac{1}{q}}.
	\end{aligned}
	\end{equation}
	Combining \eqref{3.2}, \eqref{3.3},  we get
	$$\begin{aligned}
	\mathcal{M}^p(x)	\leq&\frac{1}{8\pi}\int_{\mathbb{R}}\bigg\{
	\bigg(
	\int_{\mathbb{R}^2}|F_1(u)|^p\rho_1(u)|F_2(v)|^p\rho_2(v)G(x,u,v)du dv\bigg)\times\\
	&\times\bigg(\int_{\mathbb{R}^2} \big|(F_1\rho_1)(u)\big|\big|(F_2\rho_2)(v)\big|G(x,u,v)\mathrm{d}u\mathrm{d}v\bigg)^{\frac{p}{q}+1-p}\bigg\}dx.
	\end{aligned}$$
	Since $1\textfractionsolidus p + 1\textfractionsolidus q =1$, we infer that $\frac{p}{q}+1-p=0$.  By the assumption $\forall F_i\in L_p(\mathbb{R},\rho_i), i=\{1,2\},$ using Fubini's theorem for the right-hand side of the above equality combining with \eqref{1.8},  we deduce that	
	$$\begin{aligned}
	&\mathcal{M}^p(x)\leq\frac{1}{8\pi}
	\int_{\mathbb{R}^3}|F_1(u)|^p\rho_1(u)|F_2(v)|^p\rho_2(v)G(x,u,v)du dv dx
	\\
	&=\frac{1}{8\pi}\bigg(
	\int_{\mathbb{R}}|F_1(u)|^p\rho_1(u)du\bigg)\bigg(	\int_{\mathbb{R}}|F_2(v)|^p\rho_1(v)dv\bigg)\bigg(	\int_{\mathbb{R}}G(x,u,v)dx\bigg)
	=	\frac{1}{8\pi} 4\sqrt{2\pi}\|F_1\|_{L_p(\mathbb{R},\rho_1)} \|F_2\|_{L_q(\mathbb{R},\rho_2)}.
	\end{aligned}$$
	The proof is concluded.
\end{proof}
\noindent In case one of the  functions $\rho_1(x)$, $\rho_2(x)$ is homogenous $1$, for instance function $\rho_1(x)\equiv 1$ for all $ x\in\mathbb{R}$, and $\rho_2$ is a positive function belonging to $L_1(\mathbb{R})$ space.  By using Def.\eqref{1.5} and Gaussian integral formula \eqref{1.8},  we can compute as
$$\begin{aligned}
\big|\big(1\underset{T_c}{\overset{\gamma}{*}} \rho_2\big) (x)\big|&\leq
\frac{1}{8\pi}\int_{\mathbb{R}^2}\rho_2(v)G(x,u,v)du dv	\\	&=\frac{1}{8\pi}\bigg(\int_{\mathbb{R}}\rho_2(v)dv\bigg)\bigg( \int_{\mathbb{R}}G(x,u,v)du\bigg)=
\frac{1}{8\pi} 4\sqrt{2\pi}\bigg(\int_{\mathbb{R}}\rho_2(v)dv\bigg)=\frac{1}{\sqrt{2\pi}}\|\rho_2\|_{L_1(\mathbb{R})}<\infty.
\end{aligned}$$This shows that $(1\underset{T_c}{\overset{\gamma}{*}}\rho_2)(x)$ is well-defined $\forall x\in \mathbb{R}$. Furthermore, this operator is a bounded function $\forall x\in \mathbb{R}$ and the following estimates hold
$\|\big(1\underset{T_c}{\overset{\gamma}{*}}\rho_2\big)(x)\|_{L_{\infty}(\mathbb{R})}\leq \frac{1}{\sqrt{2\pi}}\|\rho_2\|_{L_1(\mathbb{R})}, \forall \rho_2\in L_1(\mathbb{R}),$ and
\begin{equation}\label{3.4}
\bigg|	\big(1\underset{T_c}{\overset{\gamma}{*}}\rho_2\big)(x)\bigg|^{1-\frac{1}{p}}\leq \bigg(\frac{1}{\sqrt{2\pi}}\bigg)^{1-\frac{1}{p}}\|\rho_2\|_{L_1(\mathbb{R})}^{1-\frac{1}{p}}, \forall \rho_2\in L_1(\mathbb{R}).
\end{equation}
Combining Theorem \ref{Định lý 3.1} and \eqref{3.4}, we arrive at the following corollary.
\begin{corollary}\label{Hệ quả 3.1}
	Suppose that $\rho_1\equiv 1, \forall x\in \mathbb{R}$, and $0<\rho_2\in L_1(\mathbb{R}).$ Then, for any functions $F_1\in L_p(\mathbb{R}), F_2\in L_p(\mathbb{R},\rho_2) $
	with $p > 1$, we have the following estimate
	\begin{equation}\label{3.5}
	\left\|\big(F_1 \underset{T_c}{\overset{\gamma}{*}}(F_2\rho_2)\big)\right\|_{L_p(\mathbb{R})}\leq \bigg(\frac{1}{\sqrt{2\pi}}\bigg)^{1-1\textfractionsolidus p} \|\rho_2\|_{L_1(\mathbb{R})}^{1-1\textfractionsolidus p}\|F_1\|_{L_p(\mathbb{R})}
	\|F_2\|_{L_p(\mathbb{R},\rho_2)}.
	\end{equation}
\end{corollary}
\noindent Based on \eqref{3.4},  we realize that the structure of $(\cdot \underset{T_c}{\overset{\gamma}{*}}\cdot)$ is a bounded operator in  $L_p(\mathbb{R})$ space, for any $p>1$.  Obviously,  the estimates \eqref{3.1}, and \eqref{3.4} are still valid on $L_2(\mathbb{R})$.
\section{$L_1$-solvability of Fredholm integral equation of the second kind}\label{section4}
This section will be devoted to a class of the Fredholm second kind of convolution integral equations related to $T_c$  transform \eqref{1.4} and its convolution. We will establish conditions that will
guarantee the existence and uniqueness of solutions in a closed form for these equations. We deal with a class of Fredholm equations of second type as follows (see \cite{17})
\begin{equation}\label{4.1}
f(x)-\lambda \int_{\mathcal{S}} \mathcal{K}(x,u)f(u)du=\varphi (x), \end{equation}
Here the right-hand side $\varphi (x)$ and the kernel $\mathcal{K}(x,u)$ are some known functions, $\lambda$ is a known real (complex) parameter, $\mathcal{S}$ is a piecewise-smooth surface (or line), and $f$ is an unknown function.
\begin{proposition}
	The convolution $\underset{T_c}{*}$ of two functions $f,g$ for $(T_c)$ transform is defined by
	\begin{equation}\label{4.2}
	\big(f \underset{T_c}{*}g\big)(x):=\frac{1}{2\sqrt{2\pi}}\int_{\mathbb{R}} [f(x+u)+f(x-u)]g(u)dx.
	\end{equation}
	If $f,g\in L_1(\mathbb{R})$, then $(f \underset{T_c}{*}g)$ belongs to $L_1(\mathbb{R})$ and the following factorization equality  is valid
	\begin{equation}
	\big(T_cf \big)(y)\big(T_cg \big)(y)=T_c\big(f \underset{T_c}{*}g\big)(y),
	\label{4.3}
	\end{equation}
	where $(T_c)$ transform is determined by \eqref{1.4}.\end{proposition}
\begin{proof}
	From \eqref{4.2}, for any functions $f,g \in L_1 (\R)$, we derive
	$$\begin{aligned}
	\int_{\mathbb{R}}\big|\big(f \underset{T_c}{*}g\big)(x)\big| dx&\leq \frac{1}{2\sqrt{2\pi}} \int_{\mathbb{R}^2}\big|f(x+u)+f(x-u)\big| \big|g(u)\big|du dx\\
	&\leq\frac{1}{2\sqrt{2\pi}}\bigg\{
	\int_{\mathbb{R}^2}\big|f(x+u)\big| \big|g(u)\big|du dx+
	\int_{\mathbb{R}^2}\big|f(x-u)\big| \big|g(u)\big|du dx
	\bigg\}\\
	&=\frac{1}{\sqrt{2\pi}}
	\int_{\mathbb{R}^2}\big|f(t)\big| \big|g(u)\big|du dt=\frac{1}{\sqrt{2\pi}}\|f\|_{L_1(\mathbb{R})}\|g\|_{L_1(\mathbb{R})}<\infty.
	\end{aligned}$$
	Thus, we deduce that $(f \underset{T_c}{*}g)$ belongs to $L_1(\mathbb{R}).$ Thanks to formula \eqref{1.4},  the desired relation can be achieved as follows
	$$\begin{aligned}
	(T_cf)(y)(T_cg)(y)&=\frac{1}{2\pi}\int_{\mathbb{R}^2}f(u)g(v)\cos (xu)\cos (xv)du dv
	\\&=\frac{1}{4\pi}\int_{\mathbb{R}^2}f(u)g(v)[\cos x(u+v)\cos x(u-v)]du dv
	\\&=\frac{1}{4\pi}\int_{\mathbb{R}^2}[f(\xi+y)+f(\xi-y)]g(y)\cos(\xi x)d\xi dy\\
	&=\frac{1}{\sqrt{2\pi}}\int_{\mathbb{R}}\big(f \underset{T_c}{*}g\big)(\xi)\cos(\xi x)d\xi=T_c\big(f \underset{T_c}{*}g\big)(x).
	\end{aligned}$$
\end{proof}
\noindent Our idea is to reduce the original equation \eqref{4.1} to the linear equation by using convolution $\underset{T_c}{*}$ and investigate its solvability under restriction in a domain is a real-line. We will obtain the $L_1$-solution via the simultaneous help of the factorization properties and Winner-Lévy's Theorem \cite{12}. Namely,  for kernel $\mathcal{K}(x,u)=\frac{1}{2\sqrt{2\pi}}[g(x+u)+g(x-u)]$,  choosing  $\varphi(x)=\big(f\underset{T_c}{\overset{\gamma}{*}}\psi\big)(x)$ with parameter $\lambda=-1$, and considering on infinite range $\mathcal{S} \equiv \R$. Taking into account the symmetric properties of the convolution kernel \eqref{4.2} of $(T_c)$ transform, then Eq. \eqref{4.1} can be rewritten as a convolutional equation.
\begin{equation}\label{4.4}
f(x)+\big(g\underset{T_c}{*}f\big)(x)=\big(g\underset{T_c}{\overset{\gamma}{*}}\psi\big)(x).
\end{equation}
\begin{theorem}\label{Định lý 4.1.}
	Suppose that $g,\psi$ are $L_1$-Lebesgue integrable functions over $\R$ such that $1+(T_cg)(y)\neq0$ for any $y\in\mathbb{R}$. Then, Eq.\eqref{4.4} has the unique solution in $L_1(\mathbb{R})$,  which can be represented in the form $f(x)=\big(l\underset{T_c}{\overset{\gamma}{*}}\psi\big)(x)$.  Furthermore,  the $L_1$-boundedness of the solution is ensured by the following estimate
	\begin{equation}
	\|f\|_{L_1(\mathbb{R})}\leq \frac{1}{\sqrt{2\pi}}\|l\|_{L_1(\mathbb{R})}\|\psi\|_{L_1(\mathbb{R})},
	\label{4.5}
	\end{equation}
	where $l\in L_1(\mathbb{R})$ is defined by $(T_cl)(y)=\frac{(T_cg)(y)}{1+(T_cg)(y)}, \forall y\in\mathbb{R}$ and convolution $\big(\cdot\underset{T_c}{\overset{\gamma}{*}}\cdot\big)$ is defined by \eqref{1.5}.
\end{theorem}
\begin{proof}
	Applying the $(T_c)$ transform to both sides of Eq.\eqref{4.4}, we obtain $$(T_cf)(y)+(T_cg)(y)(T_cf)(y)=\gamma(y)(T_cg)(y)(T_c\psi)(y).$$ By using \eqref{1.7},  \eqref{4.3} and under the assumption $1+(T_cg)(y)\neq0,y\in\mathbb{R},$ we deduce that
	$$(T_cf)(y)=\frac{(T_cg)(y)}{1+(T_cg)(y)}\gamma(y)(T_c\psi)(y), \forall y\in\mathbb{R}.$$
	Applying the Wiener--L\'evy's theorem \cite{12} for the $(T_c)$ transform, we can conclude that there exists a function $l\in L_1(\mathbb{R})$ such that 
	$(T_cl)(y)=\frac{(T_cg)(y)}{1+(T_cg)(y)}$, for any $y\in\mathbb{R}$. Based on the factorization identity \eqref{1.7}, this leads to $(T_cf)(y)=(T_cl)\gamma(y)(T_c\psi)(y)$ equivalent to $$(T_cf)(y)=T_c\big(l\underset{T_c}{\overset{\gamma}{*}}\psi\big)(y).$$
	Thus, $f(x)=\big(l\underset{T_c}{\overset{\gamma}{*}}\psi\big)(x)\in L_1(\mathbb{R})$ almost everywhere on $\mathbb{R}$.  Furthermore,   since $l,  \psi$ are functions belonging to the $L_1(\mathbb{R})$,  it follows from \cite{6} that $f\in L_1(\mathbb{R})$. 
	According to \eqref{1.9},  we deduce the estimate \eqref{4.5} as desired by the theorem.
\end{proof}
In addition, estimates \eqref{2.13}, \eqref{2.15a}, \eqref{2.15b}  allow us to recognize that
\begin{remark}\label{chuy2}
	\textup{	Let $p,q,r,s>1$ such that $\frac{1}{p}+\frac{1}{q}+\frac{1}{r}=2+\frac{1}{s}$ and given functions $l\in L_p(\mathbb{R}),\psi \in L_q(\mathbb{R})$,  $\gamma(y)=e^{-\frac{1}{2}y^2}\in L_r(\mathbb{R})$, which satisfy $f\in L_s(\mathbb{R})$. The upper bound product of the solution can be simplified as
		\begin{equation}\label{4.6}
		\|f\|_{L_s(\mathbb{R})} \leq\frac{1}{2\pi}\left(\sqrt{\frac{2\pi}{r}}\right)^{1\textfractionsolidus r}\|l\|_{L_p(\mathbb{R})} \|\psi\|_{L_q(\mathbb{R})},s>1.
		\end{equation}
}\end{remark} 
\begin{remark}\label{chuy3}
	\textup{i)}\textup{	Let $\psi,  \gamma \in L_1(\mathbb{R})$ and $l\in L_p(\mathbb{R})$.  Assuming that $f\in L_p(\mathbb{R})$, we deduce the boundedness of the solution as follows \begin{equation}\label{iRemark42}
		\|f\|_{L_p(\mathbb{R})} \leq \frac{1}{\sqrt{2\pi}}\|l\|_{L_p(\mathbb{R})} \|\psi\|_{L_1(\mathbb{R})}, \ \text{with}\ p>1, \gamma(y)=e^{-\frac{1}{2}y^2}.
		\end{equation}
		\textup{ii)} In case $l, \psi \in L_1(\mathbb{R})$ and $\gamma \in L_r(\mathbb{R})$.  If $f\in L_r(\mathbb{R})$ then,  we get the following estimate
		\begin{equation}\label{iiRemark42}
		\|f\|_{L_r(\mathbb{R})} \leq \frac{1}{2\pi}\left(\sqrt{\frac{2\pi}{r}}\right)^{1\textfractionsolidus r}\|l\|_{L_1(\mathbb{R})} \|\psi\|_{L_1(\mathbb{R})},\ \text{with}\ r>1, \gamma(y)=e^{-\frac{1}{2}y^2}.
		\end{equation}}
\end{remark} 
\noindent Assume that  $\psi (x)=\psi_1\rho(x)$,  where $\rho (x)\in L_1(\mathbb{R}),\psi_1\in L_1(\mathbb{R},\rho)\cap L_p(\mathbb{R},\rho)$ and $l\in L_1(\mathbb{R})\cap L_p(\mathbb{R})$ with $p>1$.  Then,  using \eqref{3.4},  we infer directly that
\begin{remark}\label{chuy4}
	\textup{For given functions $l\in L_1(\mathbb{R})\cap L_p(\mathbb{R}),0<\rho\in L_1(\mathbb{R}), \psi_1\in L_1(\mathbb{R},\rho)\cap L_p(\mathbb{R},\rho)$ with $p>1$ satisfying $\big(1\underset{T_c}{\overset{\gamma}{*}}\rho\big)$ which is well-defined.  If $f\in L_p(\mathbb{R})$, then  
		\begin{equation}\label{4.7}
		\|f\|_{L_p(\mathbb{R})}\leq \bigg(\frac{1}{\sqrt{2\pi}}\bigg)^{1-\frac{1}{p}}\|\rho\|_{L_1(\mathbb{R})}^{1-\frac{1}{p}}\|l\|_{L_p(\mathbb{R})}\|\psi\|_{L_p(\mathbb{R},\rho)}.
		\end{equation}	}
\end{remark}  
\noindent Finally, we provide an example and a thorough analysis intended to illustrate the results obtained (Theorem \ref{Định lý 4.1.} and Remarks \ref{chuy2}, \ref{chuy3} and \ref{chuy4}) to ensure their validity and applicability.

\begin{example} \rm{Let the functions $$g(x)=\sqrt{2\pi}\left\{
		\begin{aligned}
		& e^{-x} &\text {if } \,\, x>0\\
		&0 &\text {if } \, x\leq0
		\end{aligned}
		\right. \,\,\text{and} \ \,   \psi(x)=\left\{
		\begin{aligned}
		&e^{-2x} &\text {if } \, x>0\\
		&	0 & \text {if }\, x\leq0
		\end{aligned}
		\right..$$  It is easy to check that $g,\psi$ are functions belonging to $L_1(\mathbb{R}),$ and $(T_cg)(y)=\frac{1}{1+y^2}$.  
		We infer directly $1+(T_cg)(y)=1+\frac{1}{1+y^2} \ne 0,\forall y\in\mathbb{R}$.  On the other hand $\frac{(T_cg)(y)}{1+(T_cg)(y)}=\frac{1}{2+y^2}$.  Therefore, $(T_cl)(y)=\frac{1}{2+y^2}$.  By virtue of \cite{19}, we derive $$l(x)=\frac{\sqrt \pi}{2}\left\{
		\begin{aligned}
		&e^{-\sqrt{2}x} \, &\text {if}\,\, x>0\\
		&0  \,\,&\text {if } \,x\leq0
		\end{aligned}
		\right.,$$ obviously that $l$ belongs to $L_1 (\R)$ and $f(x)=\big(l\underset{T_c}{\overset{\gamma}{*}}\psi\big)(x),$ where $\gamma(y)=e^{-\frac{1}{2}y^2}$ and functions $l,\psi$ are given as above.  Invoking \eqref{4.5},  we obtain\\
		$\|f\|_{L_1(\mathbb{R})}\leq \frac{1}{\sqrt{2\pi}}\|l\|_{L_1(\mathbb{R})}\|\psi\|_{L_1(\mathbb{R})}=\frac{1}{\sqrt{2\pi}}\big\|
		\frac{\sqrt \pi}{2}e^{-\sqrt{2}x}\big\|_{L_1(\mathbb{R}_+)}\left\|
		e^{-2x}\right\|_{L_1(\mathbb{R}_+)}=\frac{1}{8}.
		$
		Assume that $p,q,r,s>1$ and satisfy the condition $\frac{1}{p}+\frac{1}{q}+\frac{1}{r}=2+\frac{1}{s}$,  functions $\psi,l$ are defined as above.   Indeed,  we obtain $l\in L_p(\mathbb{R})$, $\psi\in L_q(\mathbb{R})$ and $\gamma\in L_r(\mathbb{R})$.  Let $f\in L_s(\mathbb{R})$, based on \eqref{4.6},  we can write\\
		$\|f\|_{L_s(\mathbb{R})}\leq\frac{1}{2\pi}\left(\sqrt{\frac{2\pi}{r}}\right)^{1\textfractionsolidus r}\big\|
		\frac{\sqrt\pi}{2}e^{-\sqrt{2}x}\big\|_{L_p(\mathbb{R}_+)}\left\|
		e^{-2x}\right\|_{L_q(\mathbb{R}_+)}=\frac{1}{4\sqrt{\pi}}\left(\sqrt{\frac{2\pi}{r}}\right)^{1\textfractionsolidus r}\biggl(\frac{1}{\sqrt{2}p}\biggr)^{\frac{1}{p}}\biggl(\frac{1}{2q}\biggr)^{\frac{1}{q}}.
		$
		If $l\in L_p(\mathbb{R})$ and $\psi, \gamma\in L_1(\mathbb{R})$.  Due to the inequality \eqref{iRemark42}, it is easy to see that\\
		$\|f\|_{L_p(\mathbb{R})}\leq\frac{1}{\sqrt{2\pi}}\big\|
		\frac{\sqrt\pi}{2}e^{-\sqrt{2}x}\big\|_{L_s(\mathbb{R}_+)}\left\|
		e^{-2x}\right\|_{L_1(\mathbb{R}_+)}=\left(\frac{1}{4\sqrt{2}}\right)\left(\frac{1}{\sqrt{2}p}\right)^{\frac{1}{p}}, \forall f\in L_p(\mathbb{R}), p>1.
		$
		Since $l, \psi\in L_1(\mathbb{R})$ when $\gamma \in L_r(\mathbb{R})$ then by applying \eqref{iiRemark42},  
		we have the boundedness of the solution as follows\\
		$\|f\|_{L_r(\mathbb{R})}\leq\frac{1}{2\pi}\left(\sqrt{\frac{2\pi}{r}}\right)^{1\textfractionsolidus r}\big\|
		\frac{\sqrt\pi}{2}e^{-\sqrt{2}x}\big\|_{L_1(\mathbb{R})}\left\|
		e^{-2x}\right\|_{L_1(\mathbb{R})}=\frac{1}{8\sqrt{2\pi}}\left(\sqrt{\frac{2\pi}{r}}\right)^{1\textfractionsolidus r}, \forall f \in L_r(R), r>1.
		$\\
		We now choose functions $$\psi_1(x)=\rho(x)=\left\{
		\begin{aligned}
		&e^{-x} & \text {if} \,\,x>0\\
		&	0 & \text {if}\,\, x\leq0
		\end{aligned}
		\right..$$ It follows that $0<\rho\in L_1(\mathbb{R}),\,\psi_1\in L_1(\mathbb{R},\rho)\cap L_p(\mathbb{R},\rho),$ and $l\in L_p(\mathbb{R})$. This indicates  $$\big(1 \overset{\gamma}{\underset{T_c}{*}}\rho\big)(x)=\left\{
		\begin{aligned}
		&\big(1\overset{\gamma}{\underset{T_c}{*}}e^{-t}\big)(x) & \text {if } x>0\\
		&	0 & \text {if } x\leq0
		\end{aligned}\right..$$ With the aid of  \eqref{1.9},  we can write $\| 1\overset{\gamma}{\underset{T_c}{*}}e^{-t}\|_{L_p(\mathbb{R})}\leq \frac{1}{\sqrt{2\pi}}$.  Hence,  convolution $\big(1\overset{\gamma}{\underset{T_c}{*}}\rho\big)(x)$ is well-defined.  Let $f\in L_p(\mathbb{R})$, owing to inequality \eqref{4.7}, we have the boundedness of the solution for Eq.\eqref{4.4}, with $p>1$ as follows
		$	\|f\|_{L_p(\mathbb{R})}\leq \left(\frac{1}{\sqrt{2\pi}}\right)^{1-\frac{1}{p}}\|e^{-x}\|^{1-\frac{1}{p}}_{L_1(\mathbb{R}_+)}\big\|\frac{\sqrt\pi}{2}e^{-\sqrt{2}x}\big\|_{L_p(\mathbb{R}_+)} \|e^{-x}\|_{L_p(\mathbb{R}_+,e^{-x})}=\frac{1}{2\sqrt{2}}\biggl(\frac{\sqrt{\pi}}{p(p+1)}\biggr)^\frac{1}{p}.$}
	
\end{example}
\vskip 0.3cm

\end{document}